# A CHARACTERIZATION OF THE INFINITELY DIVISIBLE SQUARED GAUSSIAN PROCESSES

By Nathalie Eisenbaum and Haya Kaspi

*Université Paris VI–CNRS and Technion*

We show that, up to multiplication by constants, a Gaussian process has an infinitely divisible square if and only if its covariance is the Green function of a transient Markov process.

**1. Introduction.** The question of the infinite divisibility of squared Gaussian vectors is an old problem which was first raised by Paul Lévy in 1948 [10]. Given a centered Gaussian vector $(\phi_1, \ldots, \phi_p)$, when can the vector $(\phi_1^2, \ldots, \phi_p^2)$ be written as a sum of $n$ i.i.d. $p$ vectors for every $n \in \mathbb{N}$? Many authors have worked on this problem. We refer the readers to [6, 8, 13, 14, 17] and the references therein for more on this problem. In 1984, Griffiths [9] established a characterization of the $p$-dimensional centered Gaussian vectors with an infinitely divisible square. His criterion is difficult to use since it requires the computation of the signs of all the cofactors of the covariance matrix. Indeed, except for the Brownian motion (and, more generally, Gaussian Markov processes), there were no examples, in the literature, of processes satisfying this remarkable property nor examples of processes lacking it.

We have recently shown [7] that the family of fractional Brownian motions provide examples of both kinds. A fractional Brownian motion is a real-valued centered Gaussian process with a covariance given by

$$g(x,y) = |x|^\beta + |y|^\beta - |x-y|^\beta,$$

where the index $\beta$ is in $(0, 2)$. We proved that when $\beta$ is in $(0, 1]$, then the square of this process is infinitely divisible, and when $\beta$ is in $(1, 2)$, it is not. The critical index 1 corresponds to the Brownian motion. In both cases we have used the criterion of Griffiths. We have shown in [7] that if the Green function of a transient Markov process is symmetric, then it is the covariance of a Gaussian process with an infinitely divisible square. In particular, when the index $\beta$ is in $(0,1]$, the covariance of the corresponding fractional Brownian motion can be interpreted as the Green function of the symmetric stable process with index $(\beta + 1)$, killed at its hitting time of 0. To treat the second case, we have shown directly that the condition of Griffiths is not satisfied.

In view of these examples, a natural question arises. Is the representation of the covariance function of a centered Gaussian process, as the Green function of a symmetric Markov process, a necessary condition for the infinite divisibility of its square? We show in this paper that, up to a multiplication by a constant function, the answer is affirmative. The result is presented in Section 3 in the form of a necessary and sufficient condition. The proof is based this time on a criterion for infinite divisibility established by Bapat [1], which we recall in Section 2. Although, at first sight, the verification of this criterion is as difficult as that of the equivalent criterion of Griffiths (one has to check here, too, the sign of each cofactor of the covariance matrix), it allows to significantly shorten the arguments.

Gaussian processes with a covariance equal to a Green function play an important role in the study of Markov processes. Indeed, the Isomorphism theorem of Dynkin [5] provides an indentity in law connecting each of these Gaussian processes to the local time process of the corresponding symmetric Markov process. This identity has been exploited to study properties of the local time process of the Markov process, using similar properties of the Gaussian processes. We choose to mention here only two references, but many more papers on the subject are quoted in [2], for example. In [11] Marcus and Rosen have studied sample path properties of the local time process using similar properties of the Gaussian process. In [2] Bass, Eisenbaum and Shi used the Isomorphism theorem to establish the transience of the most visited sites of a symmetric stable process. The question of characterizing the Gaussian processes with a covariance that is equal to a Green function of a symmetric Markov process has been open since the Isomorphism theorem of Dynkin was first proved. The importance of an answer to this question is twofold; it will give a powerful tool for the study of these Gaussian processes, as well as that of their associated Markov processes. The results of Section 3 provide an answer to this question.

In Section 5 we show that the Brownian sheet does not have an infinitely divisible square. Recalling that linear combinations of the squared components of centered Gaussian vectors are infinitely divisible, we find this example and that of the fractional Brownian motion, with index in $(1,2)$, somewhat counterintuitive.

**2. The criterion of Bapat.** Let $A = (A_{ij})_{1 \leq i,j \leq p}$ be a $p \times p$ matrix. We write $A \geq 0$ if $A_{ij} \geq 0$ for all $i,j$.

DEFINITION 2.1. A matrix $A = (A_{ij})_{1 \leq i,j \leq p}$ is said to be an $M$-matrix if we have the following:

(i) $A_{ij} \leq 0$ for all $i \neq j$,
(ii) $A$ is nonsingular and $A^{-1} \geq 0$.

We refer the reader to [3] for the theory of $M$-matrices.

DEFINITION 2.2. A diagonal matrix is called a signature matrix if all its entries are either 1 or $-1$.

The Laplace transform, $\psi(t_1, \ldots, t_p)$ of the square of a $p$-dimensional Gaussian vector is given by

$$\psi(t) = [\det(I + GT)]^{-1/2},$$

where $t = (t_1, \ldots, t_p)$, $G$ is a positive definite $p \times p$-matrix, $T$ is the diagonal matrix with entries $T_{ii} = t_i$, and $I$ is the $p \times p$-identity matrix.

Bapat [1] obtained the following characterization of the matrices $G$ for which the Laplace transform $\psi$ is infinitely divisible, that is, for which $\psi^\delta(t)$ is a Laplace transform for any $\delta > 0$.

THEOREM A. *The Laplace transform $\psi$ is infinitely divisible if, and only if, there exists a signature matrix $S$ such that $SG^{-1}S$ is an $M$-matrix.*

REMARK 2.3. It is elementary to check, using Theorem A, that a centered three-dimensional Gaussian vector $(\eta_1, \eta_2, \eta_3)$ such that

$$\mathbb{E}(\eta_1 \eta_2) \mathbb{E}(\eta_2 \eta_3) \mathbb{E}(\eta_3 \eta_1) < 0$$

cannot have an infinitely divisible square.

**3. A necessary and sufficient condition for infinite divisibility.** In a previous work [7] on squared Gaussian processes, we have established the following result.

THEOREM B. *Let $g$ be the Green function of a strongly symmetric transient Borel right Markov process with state space $(E, \mathcal{E})$. Let $\eta$ be a centered Gaussian process with covariance $g$. Then the process $\eta^2$ is infinitely divisible.*

The set $E$ (by the assumption of a right process) is a Borel subset of a compact metric space and $\mathcal{E}$ is the $\sigma$-field of its Borel sets. Note that if $(\eta_x^2(x), x \in E)$ is an infinitely divisible squared Gaussian process, then, for any $\mathcal{E}$-measurable real valued function $d$, the process $(d^2(x)\eta_x^2(x), x \in E)$ is also infinitely divisible. We therefore have the following:

COROLLARY 3.1. *Let $g$ be the Green function of a symmetric transient Markov process on a state space $(E, \mathcal{E})$. Then for any $\mathcal{E}$-measurable real valued function $d$ on $E$, there exists a centered Gaussian process $\eta$ with covariance equal to $(d(x)g(x,y)d(y), (x,y) \in E \times E)$. The process $\eta^2$ is infinitely divisible.*

We shall first treat the case when $E$ is a finite set. Theorem 3.2 completes Corollary 3.1, for this case, by showing that its sufficient condition on the covariance is also necessary for the infinite divisibility of $\eta^2$.

THEOREM 3.2. *Let $\eta$ be a centered Gaussian vector, indexed by a finite set $E$, with a positive definite covariance function $(G(x,y), (x,y) \in E \times E)$. The vector $\eta^2$ is infinitely divisible if, and only if, there exists a real valued function $d$ on $E$ such that, for any $x, y \in E$,*

$$G(x,y) = d(x)g(x,y)\,d(y),$$

*where the function $g$ is the Green function of a transient symmetric Markov process.*

REMARK 3.3. Note that the Green function $g$ of a symmetric Markov process $X$ is always positive definite. Indeed, it is semi-positive definite (see, e.g., [11] or [7]) and it has been shown in [7], Section II, that, for any $x_1, x_2, \ldots, x_n$ in the state space of $X$, the matrix $(g(x_i, x_j))_{1 \leq i,j \leq n}$ is invertible.

Let $\eta$ be a centered Gaussian process indexed by an infinite set $E$. Then it has an infinitely divisible square if, for every finite subset $F$ of $E$, the covariance of $(\eta_x, x \in F)$ satisfies the condition of Theorem 3.2. This does not guarantee, a priori, that, for some deterministic function $d$, the covariance of $(d(x)\eta_x, x \in E)$ is the Green function of a transient symmetric Markov process. Restricting our attention to $E = \mathbb{R}$, and under the additional continuity assumption on the covariance function, the following theorem establishes the necessity of that condition.

THEOREM 3.4. *Let $\eta$ be a centered Gaussian process, indexed by $\mathbb{R}$, with a positive definite covariance function $(G(x,y), (x,y) \in \mathbb{R}^2)$. Assume that $G$* is jointly continuous. If the process $\eta^2$ is infinitely divisible, then there exists a measurable real valued positive function $d$ on $\mathbb{R}$ such that, for any $x, y \in \mathbb{R}$,

$$G(x,y) = d(x)g(x,y)\,d(y),$$

where the function $g$ is the Green function of a symmetric transient Markov process.

Theorem 3.4 is proved in Section 4. Although Theorem 3.2 was a good hint to anticipate Theorem 3.4, we could not use it directly to prove it. Our argument is based on an explicit construction of the function $d$ and of the Green function $g$ of the Markov process that will be associated with the Gaussian process $\eta$. As a first step of this construction, we show that the covariance $G$ has to be positive. The proof of Theorem 3.4 remains valid, with simple changes, when the index set $\mathbb{R}$ is replaced by a separable locally compact metric space. Moreover, we will see in Section 4 that, as a by-product of the proofs of Theorems 3.2 and 3.4, we obtain the following characterization of the associated Gaussian processes, which makes use of the definition below.

DEFINITION 3.5. A $p \times p$ matrix $A$ is said to be diagonally dominant if, for every $i = 1, \ldots, p$, $\sum_{j=1}^{p} A_{ij} \geq 0$.

THEOREM 3.6. (i) A positive definite matrix $G$ is the Green function of a finite state space transient Markov process, if, and only if, $G^{-1}$ is a diagonally dominant $M$-matrix.

(ii) Let $E$ be a separable locally compact metric space. Let $(\eta_x, x \in E)$ be a Gaussian process with a continuous positive definite covariance $(G(x,y), (x,y) \in E^2)$. Then $G$ is the Green function of a transient Markov process on $E$, if and only if, for every $x_1, x_2, \ldots, x_p$ in $E$, the inverse of the matrix $(G(x_i, x_j))_{1 \leq i,j \leq p}$ is a diagonally dominant $M$-matrix.

## 4. Proofs.

PROOF OF THEOREM 3.2. Let $E$ be the finite set $\{x_1, \ldots, x_p\}$. Theorem A guarantees the existence of a signature matrix $S$ with diagonal $S_i, i = 1, \ldots, p$ such that the covariance matrix of $(S_i \eta_{x_i}, i = 1, \ldots, p)$ is positive. We shall therefore assume, from the onset, that the covariance $G$ is positive.

We will actually prove that the covariance $G$ satisfies

(1) $$G(x,y) = \frac{D(y)}{D(x)} g(x,y),$$

where the function $g$ is the Green function of a transient Markov process and $D$ is a strictly positive function.

This will be sufficient to establish Theorem 3.2. Indeed, if $G$ satisfies (1), then we have

$$g(x,y) = \frac{D(x)}{D(y)} G(x,y).$$

Let $U$ be the corresponding potential with density $g(x,y)$ with respect to a reference measure $\mu$. That is,

$$Uf(x) = \int g(x,y) f(y) \mu(dy).$$

We then set $m(dy) = \mu(dy)/D^2(y)$. With respect to $m$, the kernel $U$ has densities $\tilde{g}(x,y) = D(x)G(x,y)D(y)$. Consequently, $\tilde{g}$ is a symmetric Green function. Setting $d = 1/D$, one obtains the necessity of the condition of Theorem 3.2.

To prove (1), let $G$ be the covariance matrix of a $p$-dimensional centered Gaussian vector with an infinitely divisible square. By Theorem A, $G^{-1}$ is an $M$-matrix. This implies (see [3]) that

$$(2) \qquad G^{-1} = cI - B,$$

where $B \geq 0$ and $c$ is strictly greater than the absolute value of any eigenvalue of $B$. Further, by [3], Chapter 6, page 137 M36, since $G^{-1}$ is an $M$-matrix, there exists a diagonal matrix $D$ such that $D_{ii} > 0$ for all $i$ and $DG^{-1}D^{-1}$ has strictly positive row sum. This means that, for any $i$,

$$\sum_{j=1}^{p} (DG^{-1}D^{-1})_{ij} > 0.$$

According to Definition 3.5, $DG^{-1}D^{-1}$ is diagonally dominant. Set $T = D(\frac{1}{c}B)D^{-1}$. Then

$$DG^{-1}D^{-1} = D(cI - B)D^{-1} = c(I - T).$$

Note that, for any $i$,

$$\sum_{j=1}^{p} T_{ij} < 1.$$

The matrix $T$ is therefore the transition matrix of a transient Markov chain $(X_n)_{n \in \mathbb{N}}$ with state space $E = \{x_1, x_2, \ldots, x_p\}$ satisfying $T_{ij} = \mathbb{P}_{x_i}(X_1 = x_j)$, and $\sum_{j=1}^{p} \mathbb{P}_{x_i}(X_1 = x_j) = 1 - \mathbb{P}_{x_i}(X_1 = \Delta)$, where $\Delta$ denotes a cemetery point. Let $\ell_\infty^{x_j}$ denote the total number of visits of $X$ to the state $x_j$. The Green function of $X$ is defined by

$$g(x_i, x_j) = \mathbb{E}_{x_i}(\ell_\infty^{x_j}).$$

It can be computed as follows:

$$\mathbb{E}_{x_i}(\ell_\infty^{x_j}) = \mathbb{E}_{x_i}\left(\sum_{n=0}^\infty \mathbb{1}_{\{X_n = x_j\}}\right) = \sum_{n=0}^\infty T_{ij}^n = (I-T)_{ij}^{-1},$$

which is defined and is finite, since, by (2) and the discussion following it, the spectral radius of $T$ is strictly smaller than 1. Hence, we can write

$$cDGD^{-1} = g;$$

that is, for every $x$, $y$ in $E$,

$$cG(x,y) = \frac{D(y)}{D(x)} g(x,y).$$

We shall now use the well-known method to turn a finite state space Markov chain into a continuous time Markov process by subordination to a Poisson process with rate $c$ (see, e.g., [4]). By its construction, this Markov process is transient with potential density (Green function) equal to $g$ and (1) is established. □

REMARK 4.1. Let $Y$ be the Markov process with potential density $g$. Equation (1) looks as if $G$ is the potential density of an $h$-path transform of $Y$. This is really the case if $D$ is excessive for $Y$. Unfortunately, this is not true in general. Consequently, we see that the collection of covariance functions that correspond to Gaussian processes with an infinitely divisible square is somewhat richer than the set of Green functions of symmetric Markov processes. This remark remains true also when the index set $E$ of the Gaussian process is infinite.

To select covariance matrices that correspond to Green functions of Markov processes, we have the condition given by Theorem 3.6(i). Indeed, assume that $G^{-1}$ is a diagonally dominant $M$-matrix. Keeping the notation of the proof of Theorem 3.2, we can choose $D = I$ and obtain $T = \frac{1}{c}B$. For any $i$, $\sum_{j=1}^p T_{ij} \leq 1$. Since the spectral radius of $T$ is strictly smaller then 1, for at least one $i$, the above inequality must be strict. Therefore, $T$ is the transition matrix of a transient Markov chain and we can conclude that $G$ is the Green function of a transient Markov process with finite state space.

To see that the condition is also necessary, consider the Green function $g$ of a transient symmetric Markov process $X$ with a finite state space $E = \{x_1, x_2, \ldots, x_p\}$. Then the inverse of the matrix $G = (g(x_i, x_j))_{1 \leq i,j \leq p}$ is a diagonally dominant $M$-matrix. Indeed $G_{ij} = \lambda_j (I-T)_{ij}^{-1}$, where $\lambda_j$ is the expected value of the exponential sojourn times in state $j$, $T$ is the transition matrix of the Markov chain $(X(S_n))_{n \in \mathbb{N}}$, and $S_n$ is the $n$th jump time of $X$.

PROOF OF THEOREM 3.4. We first show that $G$ has to be positive. For a fixed $n$, we define the function $d_n$ on $\mathbb{R}$ by

$$d_n(x) = \frac{k}{2^n} \qquad \text{if } \frac{k}{2^n} \leq x < \frac{k+1}{2^n}.$$

Let $K$ be the compact set $[a,b]$ with $a < b$, and let $d_n(K) = \{d_n(x) : x \in K\}$. The set $d_n(K)$ is finite set. Since the process $(\eta^2_{d_n(x)}, x \in K)$ is infinitely divisible, thanks to Theorem A, there exists a signature function $s_n$ (i.e., a function taking values in $\{-1, +1\}$) on $d_n(K)$ such that, for every $x, y \in K$,

$$s_n(d_n(x)) s_n(d_n(y)) G(d_n(x), d_n(y)) \geq 0.$$

Note that

(3) $\qquad s_n(d_n(x)) s_n(d_n(y)) G(d_n(x), d_n(y)) = |G(d_n(x), d_n(y))|.$

Since $G$ is continuous, we obtain, by letting $n$ tend to $\infty$, that $\lim_n (s_n(d_n(x)) \times s_n(d_n(y)))$ exists for all $(x,y)$ for which $G(x,y) \neq 0$ and, for such $(x,y)$,

$$\left( \lim_n s_n(d_n(x)) s_n(d_n(y)) \right) G(x,y) = |G(x,y)|.$$

For $(x,y)$ in $K^2$ such that $G(x,y) = 0$, there exists a finite sequence $a_1, a_2, \ldots, a_p$ such that $G(x, a_1) G(a_1, a_2) \cdots G(a_{p-1}, a_p) G(a_p, y) \neq 0$. Indeed, for $z \in K$, set $C(z) = \{u \in \mathbb{R} : G(z, u) \neq 0\}$. For each $z$, $C(z)$ is an open set and $\bigcup_{z \in K} C(z)$ is a covering of the compact set $K$. Thus, there exists a finite subcovering $C(z_1), C(z_2), \ldots, C(z_m)$ of the set $K$. The sets of the covering are not disjoint. Let $C(z_{i_1})$ be one of these sets. $C(z_{i_1})$ is a union of disjoint intervals. If $C(z_{i_1})$ does not cover $K$, then there exists $z_{i_2}$ in $\{z_1, z_2, \ldots, z_m\}$ such that $C(z_{i_2})$ covers some of the endpoints of $C(z_{i_1})$ that are in $K$. If $C(z_{i_1}) \cup C(z_{i_2})$ do not cover $K$, there exists $z_{i_3}$ in $\{z_1, z_2, \ldots, z_m\}$ such that $C(z_{i_3})$ covers some of the endpoints of $C(z_{i_1}) \cup C(z_{i_2})$ that are in $K$, and we may continue on until we exhaust all the finite covering above. Then we just have to make use of the sequence $(z_{i_1}, z_{i_2}, \ldots, z_{i_m})$ to construct $(a_1, a_2, \ldots, a_p)$ ($p \leq m$), connecting $x$ to $y$ such that $G(x, a_1) G(a_1, a_2) \cdots G(a_{p-1}, a_p) G(a_p, y) \neq 0$. Since

$$s_n(d_n(x)) s_n^2(d_n(a_1)) s_n^2(d_n(a_2)) \cdots s_n^2(d_n(a_p)) s_n(d_n(y)) = s_n(d_n(x)) s_n(d_n(y)),$$

we obtain, using (3), the existence of $\lim_n s_n(d_n(x)) s_n(d_n(y))$ for all $x$, $y$. Set

$$H(x,y) = \lim_n s_n(d_n(x)) s_n(d_n(y)).$$

The function $H$ is symmetric and, by (3),

$$H(x,y) = \text{sign}(G(x,y)),$$

and by its definition for all $x, y, z \in K$,
$$H(x,y) = H(x,z)H(z,y).$$
Hence, there exists a signature function $S_K$ on $K$ [take $H(\cdot, z_0)$ for a fixed $z_0$ in $K$] such that, for any $x, y \in K$,
$$S_K(x)S_K(y)G(x,y) = |G(x,y)|.$$
Denote by $S_n$ the function $S_{[-n,n]}$, then repeating the above argument and letting $n$ tend to $\infty$, we finally obtain the existence of $S$ satisfying
$$S(x)S(y)G(x,y) = |G(x,y)|.$$
If $S$ is not identically equal to 1, then there exists a point $x$ at which $S(x) = 1$ and $\liminf_{y \to x} S(y) = -1$. By continuity, this means that $G(x,x)$ has to be equal to 0 [because $G(x, y_n) < 0$ for a sequence $(y_n)_{n \geq 0}$ that converges to $x$]. This is excluded since $G$ is positive definite. Consequently, $G$ is positive.

We define the measure $m$ on $\mathbb{R}$ by
$$m(dy) = \left(1 \wedge \frac{1}{\sqrt{G(y,y)}}\right) e^{-|y|} \, dy.$$
Note that the measure $m$ is finite and $\int \sqrt{G(y,y)} m(dy) < \infty$. Making use of the covariance inequality ($G(x,y) \leq \sqrt{G(x,x)G(y,y)}$ for $x,y$ in $\mathbb{R}$), and the dominated convergence theorem, we see that the function

(4) $$\chi(x) = \int G(x,y)m(dy)$$

is continuous.

We now consider the fixed compact set $K = [a,b]$ with $a < b$. For any integer $n$, we keep the definition of $d_n$ introduced at the beginning of the proof. We set
$$I_n = \left\{k \in \mathbb{N} : -n2^n + 1 \leq k \leq n2^n - 1 \text{ and } \frac{k}{2^n} \in d_n(K)\right\}.$$
Define $G_n$ on the set $\{(\frac{k}{2^n}, \frac{\ell}{2^n}) : k, \ell \in I_n\}$ by
$$G_n\left(\frac{k}{2^n}, \frac{\ell}{2^n}\right) = G\left(\frac{k}{2^n}, \frac{\ell}{2^n}\right) \int_{\ell/2^n}^{(\ell+1)/2^n} m(dy).$$
Since the process $\eta^2$ is infinitely divisible and the positive measure $m$ has a support equal to $\mathbb{R}$, $G_n^{-1}$ is an $M$-matrix ($G_n$ need not be symmetric). Hence, we can write
$$G_n^{-1} = c_n I - B_n,$$
where $B_n \geq 0$ and $c_n$ is strictly greater than the absolute value of any eigenvalue of $B_n$.

Define
$$\chi_n\left(\frac{k}{2^n}\right) = \sum_{\ell \in I_n} G_n\left(\frac{k}{2^n}, \frac{\ell}{2^n}\right)$$

and let $D_n$ be the diagonal matrix $\mathrm{diag}(\chi_n(\frac{k}{2^n}), k \in I_n)$. Then, $D_n e_n = G_n e_n$, where $e_n$ is equal to $(1, 1, \ldots, 1)$. Consequently, for any $k \in I_n$,
$$\sum_{\ell \in I_n} [D_n^{-1} G_n^{-1} D_n]_{k,\ell} = \left(\chi_n\left(\frac{k}{2^n}\right)\right)^{-1} > 0.$$

Setting $T_n = \frac{1}{c_n} D_n^{-1} B_n D_n$, we have

(5) $$D_n^{-1} G_n^{-1} D_n = c_n(I - T_n)$$

and for every $k$, $\sum_{\ell \in I_n} T_n(k, \ell) < 1$. Consequently, the matrix $T_n$ is the transition matrix of a transient Markov chain.

Rewriting (5), it follows that

(6) $$G_n = D_n O_n (D_n)^{-1},$$

with $O_n = \frac{1}{c_n}(I - T_n)^{-1}$.

Let $A$ be a square matrix of size $|I_n|$ defined by $A = (A(\frac{k}{2^n}, \frac{\ell}{2^n}))_{k, l \in I_n}$. We associate with $A$ an operator on the set of the continuous functions with compact support. We denote this operator $\mathbf{A}$ and define it as follows. Let $f$ be a continuous function on $\mathbb{R}$ with a compact support, then the function $\mathbf{A}f$ is given by
$$\mathbf{A}f(x) = \sum_{\ell \in I_n} A\left(d_n(x), \frac{\ell}{2^n}\right) f\left(\frac{\ell}{2^n}\right) \quad \forall x \in \mathbb{R},$$

with the convention that $A(\frac{k}{2^n}, \frac{\ell}{2^n}) = 0$ when $k$ is outside of $I_n$.

That way we associate with $D_n$ (resp. $D_n^{-1}$, $G_n$, $O_n$) the operator $\mathbf{D_n}$ (resp. $\mathbf{D_n^{-1}}$, $\mathbf{G_n}$, $\mathbf{O_n}$). Note that we have, for every function $f$ and every $x$ in $\mathbb{R}$,

$$\mathbf{D_n} f(x) = \chi_n(d_n(x)) f(d_n(x))$$
$$(\mathbf{D_n})^{-1} f(x) = (\chi_n(d_n(x)))^{-1} f(d_n(x))$$
$$\mathbf{G_n} f(x) = \sum_{\ell \in I_n} G\left(d_n(x), \frac{\ell}{2^n}\right) f\left(\frac{\ell}{2^n}\right) \int_{\ell/2^n}^{(\ell+1)/2^n} m(dy).$$

By (5), we know that the sequence $(\mathbf{O_n})_{n \geq 0}$ satisfies
$$\mathbf{O_n} = \mathbf{D_n^{-1}} \mathbf{G_n} \mathbf{D_n}.$$

Moreover, for each $n$, $\mathbf{O_n}$ is the Green operator of a finite state space Markov process.

Let $O$ be the operator defined on $\mathcal{C}_K$ (the continuous functions with support in $K$) by

$$Of(x) = D^{-1}\mathcal{G}Df(x), \qquad x \in K,$$

where

$$Df(x) = \chi_K(x)f(x),$$
$$D^{-1}f(x) = (\chi_K(x))^{-1}f(x)$$

and

$$\mathcal{G}f(x) = \int_K G(x,y)f(y)m(dy),$$

with

$$\chi_K(x) = \int_K G(x,y)m(dy).$$

Note that, by the continuity of $G$ and the fact that $G(x,x) > 0$ and $m$ has $\mathbb{R}$ as its support, $\chi_K$ is continuous, strictly positive on $K$ and, thus, bounded below by a strictly positive constant.

LEMMA 4.2. *For every function $f \in \mathcal{C}_K$, the sequence $(\mathbf{O_n}f)_{n \geq 0}$ converges to $Of$ uniformly on $K$.*

PROOF. For any continuous function $f$ on $\mathcal{C}_K$ and any $x \in K$, we have

$$\mathbf{O_n}f(x) = \frac{1}{\chi_n(d_n(x))} \int_K G(d_n(x), d_n(y))\chi_n(d_n(y))f(d_n(y))m(dy)$$

and

$$Of(x) = \frac{1}{\chi_K(x)} \int_K G(x,y)\chi_K(y)f(y)m(dy).$$

Note that

$$\chi_n(d_n(x)) = \int_K G(d_n(x), d_n(y))m(dy).$$

Since $G$ is uniformly continuous on compacts and $m$ is finite, the sequence $(\chi_n(d_n(x)))_{n \in \mathbb{N}}$ converges uniformly on $K$ to $\chi_K(x)$. Hence, the sequence $G(d_n(x), d_n(y))\chi_n(d_n(y))f(d_n(y))$ converges uniformly on $K \times K$ to $G(x,y) \times \chi_K(y)f(y)$. □

We would like to show that $O$ is a Green operator. For this, we shall use the following lemma which, in essence, is due to Hunt (see [12], Chapter X, page 255, or [15]).

LEMMA C. *Let $V$ be a kernel on a measurable space $(E, \mathcal{E})$ such that $V$ satisfies the complete maximum principal and the function $V1$ is bounded. There exists then a sub-Markovian resolvent $(V_p)$ such that $V_0 = V$.*

As Meyer has noted in [12], page 253, if we assume that $V$ is continuous, then it is sufficient to verify the complete maximum principle for continuous functions with compact support only. More precisely, we have to verify that, for any $a \geq 0$, any two positive continuous functions with compact support, $(f, h)$, if for all $x$ in $\{h > 0\}$

$$a + Vf(x) \geq Vh(x),$$

then the inequality remains valid for all $x$ in $E$.

LEMMA 4.3. *The kernel $O$ satisfies the complete maximum principle on $K$.*

PROOF. First note that $O$ maps continuous functions with support in $K$ to continuous functions and is therefore a continuous kernel on $K$. Let $f$ and $h$ be in $\mathcal{C}_K^+$ and $a \geq 0$ and suppose that

(7) $$a + Of(x) \geq Oh(x) \quad \text{for all } x \in \{h > 0\}.$$

Recall that $\{h > 0\}$ is contained in $K$. We need to show that (7) is satisfied by all $x \in K$. Suppose that this is not true. Then there exists a constant $b > 0$ and a point $x_0$ in $K$ so that

$$a + Of(x_0) < Oh(x_0) - b.$$

Let $\varepsilon > 0$ and $N$ be such that, for any $n > N$,

$$|\mathbf{O_n} f(x) - Of(x)| < \varepsilon \quad \forall x \in K$$

and

$$|\mathbf{O_n} h(x) - Oh(x)| < \varepsilon \quad \forall x \in K.$$

Such $N$ exists by Lemma 4.2 and the fact that $f, h \in \mathcal{C}_K^+$.

By (7), for every $n > N$ and $x$ in $\{h > 0\}$,

(8) $$a + 2\varepsilon + \mathbf{O_n} f(x) \geq \mathbf{O_n} h(x).$$

Recall now that, for each $n$, $O_n$ is a potential of a Markov process on $I_n$. It follows from [16] that $O_n$ satisfies the complete maximum principle on $I_n$ and, by its definition, $\mathbf{O_n}$ satisfies the complete maximum principle on $K$. Consequently, (8) is valid for all $x$ in $K$. In particular, we have

$$a + 2\varepsilon + \mathbf{O_n} f(x_0) \geq \mathbf{O_n} h(x_0).$$

On the other hand, we have
$$a - 2\varepsilon + \mathbf{O_n}f(x_0) < \mathbf{O_n}h(x_0) - b.$$
Choosing $\varepsilon < b/4$ leads to the desired contradiction. □

In order to simplify the rest of the proof of Theorem 3.4, we use now, instead of $O$, the operator $V_K$ defined by
$$V_K f(x) = \frac{1}{\chi_K(x)} \int_K G(x,y) f(y) m(dy).$$

First we note that $V_K 1 = 1$. Further, since $V_K f = O(f/\chi_K)$, $V_K$ satisfies the complete maximum principle on $K$, by Lemma 4.3.

Let $K_n$ be the compact set $[-n, n]$ for $n \in \mathbb{N}$, and denote by $V_{K_n}$ the corresponding operator as defined above for $K$. Further, arguing as for $\chi_K$, one can show that $\chi$ defined in (4) is strictly positive.

LEMMA 4.4. *Let $V$ be the operator defined on bounded Borel functions by*

(9) $$V f(x) = \frac{1}{\chi(x)} \int_{\mathbb{R}} G(x,y) f(y) m(dy).$$

*Then there exists a sub-Markovian resolvent $(V_p)$ such that $V_0 = V$.*

PROOF. Note that $V 1 = 1$, and that $V$ maps continuous functions with compact support to continuous functions on $\mathbb{R}$. Let $f$ and $h$ be two positive continuous functions with compact supports and $a \geq 0$ and suppose that

(10) $$a + V f(x) \geq V h(x) \quad \text{for all } x \in \{h > 0\};$$

we need to show that this is satisfied for all $x$ in $\mathbb{R}$.

We denote by $H$ a compact set that contains both the compact supports of $f$ and of $h$. There exists $n_0$ such that, for $n > n_0$, $H$ is contained in $K_n$. Hence, (10) can be written as
$$a \frac{\chi(x)}{\chi_{K_n}(x)} + V_{K_n} f(x) \geq V_{K_n} h(x) \quad \text{for all } x \in \{h > 0\}.$$

We have shown, when defining $m$, that $\chi$ is a continuous function. Moreover, for any $n \in \mathbb{N}$, the function $\chi_{K_n}$ is continuous and the sequence $(\chi_{K_n})$ is increasing and converges simply to $\chi$. Consequently, by Dini's theorem, $(\chi_{K_n})$ converges uniformly to $\chi$ and $\chi(x)$ is strictly positive. Since the sequence $(\chi_{K_n})_{n > n_0}$ is bounded below by a strictly positive constant on $H$, $\chi/\chi_{K_n}$ converges to 1 uniformly on $H$. Hence, for every $\varepsilon > 0$, there exists $N$ such that, for every $n > N$,
$$a(1+\varepsilon) + V_{K_n} f(x) \geq V_{K_n} h(x) \quad \text{for all } x \in \{h > 0\}.$$

Since the operator $V_{K_n}$ satisfies the complete maximum principle on $K_n$, the above inequality is still true for $x$ in $K_n$. We now multiply each side of the inequality by $\chi_{K_n}(x)/\chi(x)$ to obtain

$$a(1+\varepsilon)\chi_{K_n}(x)/\chi(x) + Vf(x) \geq Vh(x) \qquad \text{for all } x \in K_n.$$

Since $\chi_{K_n}(x)/\chi(x) \leq 1$, we finally get

$$a(1+\varepsilon) + Vf(x) \geq Vh(x) \qquad \text{for all } x \in K_n,$$

and letting $n$ tend to $\infty$,

$$a(1+\varepsilon) + Vf(x) \geq Vh(x) \qquad \text{for all } x \in \mathbb{R}.$$

Since this is true for any $\varepsilon > 0$, (10) is established for all $x$ in $\mathbb{R}$. $\square$

The sub-Markovian resolvent $(V_p)$ allows one to construct a semi-group $(P_t)$ with $(V_p)$ as its resolvent. It will be a Ray process if we restrict our state space to a compact set, and, in general, we may need to apply a Ray Knight compactification in order to obtain a good semigroup. With this semigroup one can construct a transient Markov process on $E = \mathbb{R}$ with 0-potential equal to $V$. On the other hand, (9) implies that $m$ is a reference measure for the Markov process with potential $V$ and that the potential density $h(x, y)$ of $V$ with respect to $m$ is equal to

$$h(x, y) = \frac{1}{\chi(x)} G(x, y).$$

Setting then $\mu(dy) = \chi(y) m(dy)$, we see that the operator $V$ admits the symmetric densities $(\frac{1}{\chi(x)} G(x,y) \frac{1}{\chi(y)}, x, y \in \mathbb{R})$ with respect to $\mu$. $\square$

REMARK 4.5. Assume that $G$ is such that, for any $p$, the inverse of the covariance matrix of $(\eta_{x_1}, \eta_{x_2}, \ldots, \eta_{x_p})$ is a diagonally dominant $M$-matrix. One can then follow the steps of the proof of Theorem 3.4, without the need to define the matrices $D_n$ (we check similarly, as in Remark 4.1, that the matrix $T_n = \frac{1}{c_n} B_n$ is a transition matrix) and conclude that there exists a symmetric Markov process with potential densities equal to $(G(x,y), (x,y) \in \mathbb{R}^2)$. This, together with Theorem 3.6(i), leads to Theorem 3.6(ii).

Finally, every centered Gaussian process with infinitely divisible square is equal to a deterministic function times a Gaussian process that is associated with a Markov process. The isomorphism theorem of Dynkin can hence be easily adapted to incorporate the deterministic functions, so that it can be used in studying any Gaussian processes with infinitely divisible square.

**5. The case of the Brownian sheet.** A Brownian sheet is a centered Gaussian process $(W_{x,s}, x \geq 0, s \geq 0)$ with a covariance given by

$$\mathbb{E}(W_{x,s} W_{y,t}) = (x \wedge y)(s \wedge t).$$

Remember that, up to a multiplicative constant, $(x \wedge y; x, y \in \mathbb{R}_+)$ is the Green function of the linear Brownian motion killed at its first hitting time of 0. Hence, this covariance is the product of two Green functions and one may ask whether $(W_{x,s}^2, x \geq 0, s \geq 0)$ is infinitely divisible. The answer is given in the following proposition.

PROPOSITION 5.1. (i) *For every $(x_i, s_i)_{1 \leq i \leq 3}$ of $\mathbb{R}_+^6$, the vector $(W_{x_i,s_i}^2, 1 \leq i \leq 3)$ is infinitely divisible.*
(ii) *The process $(W_{x,s}^2, x \geq 0, s \geq 0)$ is not infinitely divisible.*

PROOF. Using Griffiths criterion (or Bapat's criterion), it is easy to check that a sufficient condition for the infinite divisibility of the square of a three-dimensional Gaussian vector, indexed by $\{1, 2, 3\}$ and with a covariance $g$, is

(11) $$g(i,j)g(k,k) \geq g(i,k)g(j,k),$$

for any choice of $i, j, k$ in $\{1, 2, 3\}$.

Since any Green function satisfies (11), so does the covariance function of the Brownian sheet, as a product of two Green functions.

By Bapat's criterion, we know that a covariance matrix $G$, such that $G \geq 0$, corresponds to an infinitely divisible square Gaussian vector if $G^{-1}$ is an $M$-matrix. We choose $(x_i, s_i)_{1 \leq i \leq 4}$ such that $0 < x_1 < x_3 < x_2 < x_4$ and $0 < s_4 < s_1 < s_3 < s_2$. Let $G$ be the matrix $((x_i \wedge x_j)(s_i \wedge s_j))_{1 \leq i,j \leq 4}$. We compute the coefficient $G_{1,2}^{-1}$:

$$G_{1,2}^{-1} = x_1 x_3 s_4 (x_2 - x_3)(s_3 - s_1) s_4 > 0.$$

Hence, $G^{-1}$ is not an $M$-matrix and the corresponding Gaussian vector does not have an infinitely divisible square. □

**Acknowledgments.** We would like to thank Gennady Samorodnitsky and Jay Rosen for helpful discussions and comments which helped us to improve a previous version of this work.

LABORATOIRE DE PROBABILITÉS
UNIVERSITÉ PARIS VI–CNRS
4 PLACE JUSSIEU
75252 PARIS CEDEX 05
FRANCE
E-MAIL: nae@ccr.jussieu.fr

INDUSTRIAL ENGINEERING
AND MANAGEMENT
TECHNION
HAIFA
ISRAEL 32000
E-MAIL: iehaya@tx.technion.ac.il